\magnification=1200
\input amstex
\documentstyle{amsppt}
\hsize=165truemm
\vsize=227truemm
\def\C{{{\Bbb C}}}
\def\Hilb{{{\Cal H}\kern -0.25ex{\italic ilb\/}}}
\def\p#1{{{\Bbb P}^{#1}_{\C}}}
\def\Ofa#1{{{\Cal O}_{#1}}}

\def\PGL{\operatorname{PGL}}
\def\Aut{\operatorname{Aut}}
\def\mapright#1{\mathbin{\smash{\mathop{\longrightarrow}
\limits^{#1}}}}

\topmatter
\title
On the birational geometry of moduli spaces of pointed curves
\endtitle
\author
E. Ballico, G. Casnati, and C. Fontanari
\endauthor

\address
Edoardo Ballico, Dipartimento di Matematica, Universit\`a degli Studi 
di Trento, Via Sommarive 14, 38050 Povo, Italy
\endaddress

\email
ballico\@science.unitn.it
\endemail

\address
Gianfranco Casnati, Dipartimento di Matematica, Politecnico di Torino,
c.so Duca degli Abruzzi 24, 10129 Torino, Italy
\endaddress

\email
casnati\@calvino.polito.it
\endemail

\address
C. Fontanari, Dipartimento di Matematica, Politecnico di Torino,
c.so Duca degli Abruz\-zi 24, 10129 Torino, Italy
\endaddress

\email
claudio.fontanari\@polito.it
\endemail


\thanks
This research was partially supported by MIUR and GNSAGA of INdAM 
(Italy).
\endthanks

\abstract
We prove that the moduli space ${\Cal M}_{g,n}$ of smooth curves 
of genus $g$ with $n$ marked points is rational for $g=6$ and $1 \le n \le 8$, 
and it is unirational for $g=8$ and $1 \le n \le 11$,  $g=10$ and $1 \le n \le 3$,
$g=12$ and $n = 1$.  
\endabstract

\endtopmatter

\document

\head
0. Introduction and notation
\endhead

Let ${\Cal M}_{g,n}$ be the (coarse) moduli space of smooth curves of genus $g$ 
with $n$ marked points defined over the field $\C$ of complex numbers and $\overline{\Cal M}_{g,n}$ its Deligne--Mumford compactification. Even if one is 
mainly interested in the geometry of $\overline{\Cal M}_g$ (i.e. the case $n=0$), it is 
worthwhile to consider pointed moduli spaces as well. Indeed, the boundary components 
of $\overline{\Cal M}_g$ are images via natural gluing morphisms of pointed moduli 
spaces of lower genus.

As an example of this philosophy applied in the specific context 
of birational geometry, we just mention the nice contribution [B--V]: in order to prove 
that $\overline{\Cal M}_{15}$ is rationally connected, a key ingredient (see Theorem~4.7) 
is provided by the fact that $\overline{\Cal M}_{14,2}$ is unirational. 

Several results about the Kodaira dimension $\kappa(\overline{\Cal M}_{g,n})$ of $\overline{\Cal M}_{g,n}$ are already known thanks to the papers [Be], [B--F], [B--V], [Fa] and [Lg]. In particular, in the above contributions it is proved that $\overline{\Cal M}_{g,n}$ is unirational for $1\le n\le \sigma(g)$ and $\kappa(\overline{\Cal M}_{g,n})\ge0$ for $n\ge\tau(g)$ where the values of $\sigma$ and $\tau$ are as follows
\bigbreak
\centerline{
\vbox{\offinterlineskip
\hrule
\halign{&\vrule#&\strut\quad\hfil#\quad\cr
height2pt
&\omit&
&\omit&
&\omit&
&\omit&
&\omit&
&\omit&
&\omit&
&\omit&
&\omit&
&\omit&
&\omit&
&\omit&
&\omit&
&\omit&
&\omit&
&\omit&
&\omit&
&\omit&
&\omit&
&\omit&
&\omit&
&\omit&\cr
&\hfill\hskip-3truemm$g$\hskip-3truemm\hfill&
&\hfill\hskip-3truemm$1$\hskip-3truemm\hfill&
&\hfill\hskip-3truemm$2$\hskip-3truemm\hfill&
&\hfill\hskip-3truemm$3$\hskip-3truemm\hfill&
&\hfill\hskip-3truemm$4$\hskip-3truemm\hfill&
&\hfill\hskip-3truemm$5$\hskip-3truemm\hfill&
&\hfill\hskip-3truemm$6$\hskip-3truemm\hfill&
&\hfill\hskip-3truemm$7$\hskip-3truemm\hfill&
&\hfill\hskip-3truemm$8$\hskip-3truemm\hfill&
&\hfill\hskip-3truemm$9$\hskip-3truemm\hfill&
&\hfill\hskip-3truemm$10$\hskip-3truemm\hfill&
&\hfill\hskip-3truemm$11$\hskip-3truemm\hfill&
&\hfill\hskip-3truemm$12$\hskip-3truemm\hfill&
&\hfill\hskip-3truemm$13$\hskip-3truemm\hfill&
&\hfill\hskip-3truemm$14$\hskip-3truemm\hfill&
&\hfill\hskip-3truemm$15$\hskip-3truemm\hfill&
&\hfill\hskip-3truemm$16$\hskip-3truemm\hfill&
&\hfill\hskip-3truemm$17$\hskip-3truemm\hfill&
&\hfill\hskip-3truemm$18$\hskip-3truemm\hfill&
&\hfill\hskip-3truemm$19$\hskip-3truemm\hfill&
&\hfill\hskip-3truemm$20$\hskip-3truemm\hfill&
&\hfill\hskip-3truemm$21$\hskip-3truemm\hfill&\cr
height2pt&\omit&
&\omit&&\omit&&\omit&&\omit&&\omit&&\omit&&\omit&&\omit&&\omit&&\omit&&\omit&&\omit&&\omit&&\omit&&\omit&&\omit&&\omit&&\omit&&\omit&&\omit&&\omit&\cr
\noalign{\hrule}
&\hfill\hskip-3truemm$\sigma(g)$\hskip-3truemm\hfill&
&\hfill\hskip-3truemm$10$\hskip-3truemm\hfill&
&\hfill\hskip-3truemm$11$\hskip-3truemm\hfill&
&\hfill\hskip-3truemm$14$\hskip-3truemm\hfill&
&\hfill\hskip-3truemm$15$\hskip-3truemm\hfill&
&\hfill\hskip-3truemm$12$\hskip-3truemm\hfill&
&\hfill\hskip-3truemm$15$\hskip-3truemm\hfill&
&\hfill\hskip-3truemm$11$\hskip-3truemm\hfill&
&\hfill\hskip-3truemm$5$\hskip-3truemm\hfill&
&\hfill\hskip-3truemm$8$\hskip-3truemm\hfill&
&\hfill\hskip-3truemm\hskip-3truemm\hfill&
&\hfill\hskip-3truemm$10$\hskip-3truemm\hfill&
&\hfill\hskip-3truemm\hskip-3truemm\hfill&
&\hfill\hskip-3truemm\hskip-3truemm\hfill&
&\hfill\hskip-3truemm$2$\hskip-3truemm\hfill&
&\hfill\hskip-3truemm\hskip-3truemm\hfill&
&\hfill\hskip-3truemm\hskip-3truemm\hfill&
&\hfill\hskip-3truemm\hskip-3truemm\hfill&
&\hfill\hskip-3truemm\hskip-3truemm\hfill&
&\hfill\hskip-3truemm\hskip-3truemm\hfill&
&\hfill\hskip-3truemm\hskip-3truemm\hfill&
&\hfill\hskip-3truemm\hskip-3truemm\hfill&\cr
height2pt&\omit&
&\omit&&\omit&&\omit&&\omit&&\omit&&\omit&&\omit&&\omit&&\omit&&\omit&&\omit&&\omit&&\omit&&\omit&&\omit&&\omit&&\omit&&\omit&&\omit&&\omit&&\omit&\cr
\noalign{\hrule}
&\hfill\hskip-3truemm$\tau(g)$\hskip-3truemm\hfill&
&\hfill\hskip-3truemm$11$\hskip-3truemm\hfill&
&\hfill\hskip-3truemm\hskip-3truemm\hfill&
&\hfill\hskip-3truemm\hskip-3truemm\hfill&
&\hfill\hskip-3truemm$16$\hskip-3truemm\hfill&
&\hfill\hskip-3truemm$15$\hskip-3truemm\hfill&
&\hfill\hskip-3truemm$16$\hskip-3truemm\hfill&
&\hfill\hskip-3truemm$14$\hskip-3truemm\hfill&
&\hfill\hskip-3truemm$14$\hskip-3truemm\hfill&
&\hfill\hskip-3truemm$13$\hskip-3truemm\hfill&
&\hfill\hskip-3truemm$11$\hskip-3truemm\hfill&
&\hfill\hskip-3truemm$11$\hskip-3truemm\hfill&
&\hfill\hskip-3truemm$13$\hskip-3truemm\hfill&
&\hfill\hskip-3truemm$11$\hskip-3truemm\hfill&
&\hfill\hskip-3truemm$10$\hskip-3truemm\hfill&
&\hfill\hskip-3truemm$10$\hskip-3truemm\hfill&
&\hfill\hskip-3truemm$9$\hskip-3truemm\hfill&
&\hfill\hskip-3truemm$9$\hskip-3truemm\hfill&
&\hfill\hskip-3truemm$9$\hskip-3truemm\hfill&
&\hfill\hskip-3truemm$7$\hskip-3truemm\hfill&
&\hfill\hskip-3truemm$6$\hskip-3truemm\hfill&
&\hfill\hskip-3truemm$4$\hskip-3truemm\hfill&\cr
height2pt&\omit&&\omit&&\omit&&\omit&&\omit&&\omit&&\omit&&\omit&&\omit&&\omit&&\omit&&\omit&&\omit&&\omit&&\omit&&\omit&&\omit&&\omit&&\omit&&\omit&&\omit&&\omit&\cr}
\hrule}}
\bigbreak
\noindent Notice that only for few values of $g$ (namely
$g=1,4,6$ and $11$) such results are tight.

For small $g$ even the rationality of ${\Cal M}_{g,n}$ can be proved. More precisely ${\Cal M}_{g,n}$ is rational for $g=1$ and $1\le n\le 10$ (see [Be]), $g=2$ and $1\le n\le 12$, $g=3$ and $1\le n\le 14$, $g=4$ and $1\le n\le 15$, $g=5$ and $1\le n\le 12$ (see [C--F]).

In the present paper our first main theorem concerns the rationality of ${\Cal M}_{6,n}$.

\proclaim{Theorem A}
The moduli space ${\Cal M}_{6,n}$ is rational for $1\le n\le 8$.
\endproclaim

Our key idea is to mix the techniques in [S--B], where the rationality of ${\Cal M}_6$ 
was originally proven, with those introduced in [C--F]. 

In higher genera we are able to improve the result of [Lg] in some cases, according to our second main theorem.

\proclaim{Theorem B}
The moduli space ${\Cal M}_{g,n}$ is unirational for $g=8$ and $1\le n\le 8$,
$g=10$ and $1\le n\le 3$, $g=12$ and $n=1$.
\endproclaim

Our methods are miscellaneous in nature and involve different ad hoc constructions
for each genus. In particular, for $g=8$ we embed the general curve on a quadric 
surface via two distinct $g^1_5$; for $g=10$ we apply the classical construction 
of Beniamino Segre revisited in [A--C]; finally, for $g=12$ we follow the lines 
of the unirationality arguments presented in [Ve]. 

We are grateful to Igor Dolgachev, who shared with us his deep insight into 
rationality questions in a stimulating and inspiring conversation held in Turin
during the School (and Workshop) on Cremona transformations (September 2005). 

\subhead
Notation
\endsubhead
We work over the field $\C$ of complex numbers.  Projective $r$--space will be denoted by $\p r$.

If $X$ is a projective scheme, then a divisor $D$ on $X$ is always intended
to be a Cartier divisor and $\vert D\vert$ denotes the complete linear system
of divisors linearly equivalent to $D$.

A curve $C$ is a projective scheme of dimension $1$. If $D$ is a divisor on $C$ then $\varphi_{\vert D\vert}$ denotes the map induced on $C$ by $\vert D\vert$. If $C, D \subset \p2$ are distinct plane  curves, then $C \cdot D$ denotes their intersection cycle.

Let ${\Cal M}_{g,n}$ be the coarse moduli space of smooth, projective $n$--pointed curves (briefly $n$--pointed curves in what follows) of
genus $g$. A point of ${\Cal M}_{g,n}$ is then the isomorphism class of a $(n+1)$--tuple $(C,p_1,\dots,p_n)$, where $C$ is a  smooth and connected curve of genus $g$ and $p_1,\dots,p_n\in C$.
As usual, we set ${\Cal M}_g := {\Cal M}_{g,0}$.
A general curve of genus $g$ is a curve corresponding to the general  point of the irreducible scheme ${\Cal M}_g$.

We denote isomorphisms by $\cong$ and birational equivalences by $\approx$.
For other definitions, results and notation we always refer to [Ha].

\head
1. The rationality of ${\Cal M}_{6,1}$
\endhead
Let $(C,p_1)\in{\Cal M}_{6,1}$ be a general point: in particular $C$ is a general curve of genus $6$. Following [SB], we notice that 
the canonical system $\vert K\vert$ embeds $C$ in $\p5$ as a curve $C_K$ of degree $10$ as the intersection of a unique del Pezzo surface $S_5\subseteq\p5$ of degree $5$ with a quadric.
We recall that such an $S$ is the blow up of $\p2$ at four general points embedded in $\p5$ via the linear system $\Sigma$ of cubic divisors in $\p2$ through the four points.

Since each two sets of four general points in $\p2$ are projectively equivalent, we can always assume that such four points are the fundamental points $E_1,E_2,E_3,E_4$ and let $s\colon S\to\p2$ be the corresponding blow up morphism. In particular, we can fix once and for all the surface $S$ and such a representation of $S$ allows us to consider insted of $C_K$ its plane model $\overline{C}:=s(C_K)\subseteq\p2$ as plane sextic curve with nodes at the points $E_1,E_2,E_3,E_4$.

We finally recall that $\Aut(S)\cong{\frak S}_5$: in the plane representation of $S$, it coincides with the rational representation of ${\frak S}_5$ as the group generated by the projectivities fixing the set $E:=\{\ E_1,E_2,E_3,E_4\ \}$ (which is isomorphic to ${\frak S}_4$) and by the standard quadratic transformation
$$
\mu(x_1,x_2,x_3)=(x_2x_3,x_1x_3,x_1x_3).
$$

\remark{Remark 1.1}
Let $\Psi\colon C\mapright\sim C'$ be an isomorphism of abstract $1$--pointed curves  of genus $6$ mapping the point $p_1\in C$ to $p_1'\in C'$. It carries the canonical system $K$ of $C$ onto the one $K'$ of $C'$. Thus $\Psi$ induces a projectivity of the canonical space mapping $C_K$ to $C'_{K'}$: obviously such a projectivity must map the unique del Pezzo surface through $C$, which is $S$, to the one through $C'$, which is again $S$.

Thus $\Psi$ induces via $s$ an element $\sigma\in{\frak S}_5$ mapping $\overline{C}$ to $\overline{C'}:=s(C'_{K'})$ and $A_1:=s(p_1)$ to $A_1':=s(p_1')$. 

Conversely, each $\sigma\in\Aut(S)\cong{\frak S}_5$ mapping $\overline{C}$ to $\overline{C'}$ and $A_1$ to $A_1'$ is induced by an isomorphism $\Psi\colon C\mapright\sim C'$ such that $\Psi(p_1)=p_1'$, $i=1,\dots,n$. 
\endremark
\medbreak

If $V\subseteq{\Bbb C}[x_1,x_2,x_3]_6$ is the subspace corresponding to curves having double points at the points of $E:=\{\ E_1,E_2,E_3,E_4\ \}$, the map $s$ induces an isomorphism $H^0\big(S,\Ofa S(2)\big)\cong V$ as representations of ${\frak S}_5$. Consider now the incidence variety
$$
X':=\{\ (f,A_1)\in{\Bbb P}(V)\times \p2\ \vert\ f(A_1)=0\ \},
$$
the irreducible open subscheme $X\subseteq X'$ and the open set ${\Cal U}\subseteq\p2 \setminus E$ defined in Lemma~2.5.1 of [C--F]. Since the general element of ${\Bbb P}(V)$ is an irreducible sextic having nodes at the points of $E$, then the same is true for the curve represented by general points in $X$.

The above description yields a rational map $m_{6,1}\colon X\dashrightarrow{\Cal M}_{6,1}$. Remark~1.1 implies that its fibres are exactly the orbits with respect to the above described action of ${\frak S}_5$.  The same argument used in Section~2.5 of [C--F] implies that $m_{6,1}$ is dominant.

As a consequence, there exists a birational equivalence ${\Cal M}_{6,1}\approx X/{\frak S}_5$ and we will obtain the rationality of ${\Cal M}_{6,1}$ by proving the rationality of the quotient on the right. 

We make the following two assertions.

\proclaim{Claim 1.2}
The action of ${\frak S}_5$ on $S$ is almost free.
\endproclaim 

Assuming the above Claim we are now ready to prove the following

\proclaim{Theorem A for $g=6$ and $n=1$}
The moduli space ${\Cal M}_{6,1}$ is rational.
\endproclaim
\demo{Proof}
From now on we will make the obvious identification of $\p2 \setminus E$ with an open subset of $S$. We first recall (see [SB], Corollary~3) that $V$ contains both the trivial and the signature representations of ${\frak S}_5$. In particular ${\Bbb P}(V)\times S$ contains a unisecant, whose intersection with $X$ is a unisecant itself, hence there exists a ${\frak S}_5$--equivariant birational equivalence  $X\approx S\times{\Bbb C}^{14}$.
Now, by taking into account Claim~1.2 above and Section~4 of [Do], it follows that the first projection $X\approx S\times{\Bbb C}^{14}\to S$ induces a vector bundle structure $X/{\frak S}_5\approx (S\times{\Bbb C}^{14})/{\frak S}_5\to S/{\frak S}_5$. Since $S$ is rational then $S/{\frak S}_5$ is a unirational surface, thus it is rational by a Theorem of Castelnuovo.
\qed
\enddemo

Next we turn to the Claim.

\demo{Proof of Claim~1.2}
We recall that $S$ contains exactly five pencils of conics. In our plane representation of $S$ they correspond to the four pencils of lines through each one of the points $E_i$ and to the pencil of conics with base points $E_0,E_1,E_2,E_3$. In particular, two distinct conics in the same pencil do not intersect and conics in different pencils intersect exactly at one point.

Choose a conic $D\subseteq S$ in the pencil corresponding to the plane conics through the points of $E$. Consider the closure $F$ of the set $(\bigcup_{\sigma\in{\frak S}_5}\sigma(D))\setminus D$ and notice that $F\cap D$ is a finite number of points. It follows that if $\sigma\in{\frak S}_5$ is in the stabilizer of a point $D\setminus F$ then $\sigma(D)=D$, hence $\sigma$ fixes the pencil $\vert D\vert$. 

Via $s$, the permutation $\sigma$ corresponds to an element of the subgroup of ${\frak S}_5$ generated by the projectivities fixing the set $E$, which is isomorphic to ${\frak S}_4$. Such a representation of ${\frak S}_4$ is the standard permutation representation, whence it is clear that for a general choice of $p\in D\setminus F$ the stabilizer of $p$ is trivial.
\qed
\enddemo

\head
2. The rationality of ${\Cal M}_{6,n}$ for $2\le n\le 8$
\endhead
Let $(C,p_1,\dots,p_n)\in{\Cal M}_{6,n}$ be general: then $C$ is general in ${\Cal M}_6$, thus it is neither hyperelliptic
nor trigonal.
As pointed out in the previous section for each such $C$, the canonical system $\vert K\vert$ is very ample and
we obtain an embedding $\varphi\colon C\to \widetilde{C}\subseteq\p5$
of degree $10$.

For every choice of general points $p_1,p_2\in C$ the linear system
$\left\vert K-2p_{n-1}-p_n\right\vert$ is a $g_7^2$. The induced map $\varphi\colon C\to \overline{C}\subseteq \p2$ coincides with the projection of $\widetilde{C}$ from the plane spanned by the tangent line to $C$ at the point $p_{n-1}$ and $p_n$.

Its image $\overline{C}$ is a septic in $\p2$ passing through the points
$A_i:=\varphi(p_i)$, $i=1,\dots,n$. 

\proclaim{Proposition 2.1}
Under the above hypotheses, the curve $\overline{C}$ is an integral septic with nine nodes $N_1,\dots,N_9$. Moreover there exists a unique smooth cubic $E_{\overline{C}}$ such that
$$
\overline{C}\cdot E_{\overline{C}}=2\sum_{i=1}^9N_i+2A_{n-1}+A_n.
$$
Every other integral septic $D$ with nine nodes at $N_1,\dots,N_9$ and tangent in $A_{n-1}$ to $E_{\overline{C}}$ cuts out residually on such cubic the point $A_n$.

Conversely, each integral septic $D\subseteq\p2$ with exactly nine nodes $N_1,\dots,N_9$ lying on an smooth cubic $E\subseteq\p2$ such that $D\cdot E=2\sum_{i=1}^9N_i+2A_{n-1}+A_n$ for two suitable points $A_{n-1},A_n\not\in\{\ N_1,\dots,N_9\ \}$ is obtained in the above way.
\endproclaim
\demo{Proof}
We start by proving the second part of the statement. Let $\ell\subseteq\p2$ be a line. Since the embedding $D\subseteq\p2$ is given by the linear system
$$
\vert\ell\cdot D\vert\subseteq \vert \ell\cdot D+E\cdot D-2\sum_{i=1}^9N_i-2A_{n-1}-A_n\vert\subseteq\vert K-2A_{n-1}-A_n\vert,
$$
it suffices to check that $\dim\vert K-2A_{n-1}-A_n\vert=2$. If not, $\left\vert 2A_{n-1}-A_n\right\vert$ would be a $g_3^1$ on $D$, but this is impossible
by the following classical argument (see [E--C], Libro Quinto, p. 106). For each $B_1+B_2+B_3\in \left\vert 2A_{n-1}-A_n\right\vert$, consider the linear system $\Sigma_{B_1+B_2+B_3}$ of quartics through the points $B_1,B_2,B_3,N_1,\dots,N_9$ and notice that
$$
\dim\Sigma_{B_1+B_2+B_3}\ge\dim\vert K-B_1-B_2-B_3\vert=\dim\vert K-2A_{n-1}-A_n\vert\ge3.
$$
The general element $Q\in\Sigma_{B_1+B_2+B_3}$ is smooth and $\Sigma_{B_1+B_2+B_3}$ cuts out on $Q$ residually to $B_1+B_2+B_3+N_1+\ldots+N_9$ a $g_4^r$ with $r\ge2$, which has to be the canonical  $g_4^2$ on $Q\subseteq\p2$ cut out by the lines (since $Q$ is smooth). It follows that $B_1,B_2,B_3,N_1,\dots,N_9$ lie on a cubic due to a well--known corollary of Noether's $A\Phi+B\Psi$ Theorem (see [Wa], Theorem 7.7). Thus we should have infinitely many cubics through $N_1,\dots,N_9$, one for each choice of $B_1+B_2+B_3$, a contradiction if the points $N_1,\dots,N_9$ are general.

Clearly $\overline{C}$ is singular and by the genus formula
$$
\sum_{P\in \overline{C}}m_P(m_P-1)=18,
$$
$m_P$ being the multiplicity of $P$ on the curve $\overline{C}$. Thus $\overline{C}$ carries at most four--fold points as singularities. Since $C$, being general, is not trigonal, then $\overline{C}$ carries at most triple points.

Simple parameter computations show that the families of septics with at least either one triple point or a double point with an infinitely near double point have dimension at most $25$, while the family of septic carrying nine pairwise distinct double points has dimension $26$. Since we have proved above that each plane septic with such a last configuration of singularities is the projection from a suitable secant plane of its canonical model, it follows that for a general curve $C$ of genus $6$ and for general points $p_1,p_2\in C$, the singularities of $\overline{C}$ are exactly nine pairwise distinct double points. We denote them by $N_1,\dots,N_9$.

Such nine nodes $N_1,\dots,N_9$ lie necessarily on at least one cubic. Since through nine general points in the plane there is always a septic having them as nodes (it is obtained as projection of its canonical model again by the first part of the proof), it follows that for a general $C$ such a cubic is uniquely determined: we denote it by $E_{\overline{C}}$.

The curve $E_{\overline{C}}$ is integral. If not, it is the union of a line and a conic. Due to degree reasons, the line and the conic contain at most three and seven nodes respectively, thus there are only two possible distinct configurations. Again a parameter computation shows that both such configurations yield families of dimension at most $24$ and we can argue as in the previous cases.
In the same way one can deduce that $E_{\overline{C}}$ is smooth.

Next, we check that $A_{n-1},A_n\in E_{\overline{C}}$. By construction, for each general line $\ell\subseteq \p2$ the divisor $\ell\cdot \overline{C}+2A_{n-1}-A_n$ is in the canonical system. Since the canonical system is cut out on $\overline{C}$ by the adjoints of degree $4$ we have the linear equivalence on $\overline{C}$
$$
\ell\cdot \overline{C}+2A_{n-1}-A_n\sim \ell\cdot \overline{C}+E_{\overline{C}}\cdot \overline{C}-2\sum_{i=1}^9N_i.
$$
Thus $2A_{n-1}-A_n\sim E_{\overline{C}}\cdot \overline{C}-2\sum_{i=1}^9N_i$ and since  the general curve $\overline{C}$ is not trigonal we have
$E_{\overline{C}}\cdot \overline{C}=2A_{n-1}-A_n+2\sum_{i=1}^9N_i$, i.e.
$A_{n-1},A_n\in E_{\overline{C}}$ and the curves $E_{\overline{C}}$ and $\overline{C}$ are tangent at $A_{n-1}$.

Finally, consider the linear system $\Cal{L}$ of plane septic curves passing
doubly through $N_1,\dots,N_9$ and and tangent to $E_{\overline{C}}$ at the point $A_{n-1}$: by Bezout's theorem,
every element of $\Cal{L}$ intersects the cubic $E_{\overline{C}}$ in another point
$A_n$. This point has to be fixed, otherwise $\Cal{L}$ would cut a $g^1_1$ on
$E_{\overline{C}}$.
\qed
\enddemo

Let $V_{7,9} \subseteq \vert \Cal O_\p2(7) \vert$ be the Severi
variety of plane curves of degree $7$ with exactly $9$ nodes as
singularities. We consider the closure $V_{7,9}^{tang} \subseteq V_{7,9}$ of the locus of those curves which are tangent to the cubic through their nodes at some point.

Proposition 2.1 shows that for each general $D\in V_{7,9}^{tang}$ then $D=\overline{C}$, thus we obtain two points $A_{n-1},A_n\in D\subseteq\p2$ where $A_{n-1}$ is the point of tangency of $D$ with the cubic through its nodes, and $A_n$ the residual intersection point. Hence we have a rational map
$V_{7,9}^{tang}\dashrightarrow \p2$ sending the general $D=\overline{C}$ to $A_{n-1}$: let us consider also the rational map $(\p2)^{n-1}\dashrightarrow \p2$ induced by the projection on the last factor. Thus we may consider the scheme
$$
W_{7,9}^n:=\{\ (D,A_1,\dots,A_{n-1})\ \vert\ \text{$A_1,\dots,A_{n-1}\in D$}\ \}\subseteq V_{7,9}^{tang}\times_{\p2}(\p2)^{{n-1}}.
$$
Now Proposition 2.1 yields the following.

\proclaim{Corollary 2.2}
The rational map $\psi\colon W_{7,9}^n\dashrightarrow{\Cal M}_{6,n}$ sending
$(D,A_1,\dots, A_{n-1})$ to its pointed normalization $(C,p_1,\dots,p_n)$ is dominant.

The fibres are the orbits of the natural action of $\PGL_3$, hence there
is a natural birational equivalence $W_{7,9}^n/\PGL_3\approx{\Cal M}_{6,n}$.
\endproclaim
\demo{Proof}
It remains to prove the assertion on the fibres. It is obvious that fibre over a point $(C,p_1,\dots,p_n)$ contains the $\PGL_3$--orbits of the plane pointed curve $(\overline{C},A_1,\dots,A_{n-1})$.

Conversely, let $\Psi\colon C\mapright\sim C'$ be an isomorphism of abstract $n$--pointed curves  of genus $6$ mapping the point $p_i\in C$ to $p_i'\in C'$. It carries the canonical system $K$ of $C$ onto the one $K'$ of $C'$. Thus $\Psi$ induces a projectivity of the canonical space mapping the canonical models $\widetilde{C}$ and $\widetilde{C}'$ one onto the other: obviously such a projectivity must map the tangent line at $p_{n-1}$ onto the tangent line at $p_{n-1}'$, thus $\Psi$ maps the linear system $\vert K-2p_{n-1}-p_n\vert$ to $\vert K'-2p_{n-1}'-p_n'\vert$, hence it induces via the projections $\varphi$ and $\varphi'$ a projectivity of $\p2$ sending $\overline{C}$ to $\overline{C}'$ and $A_i:=\varphi(p_i)$ to $A_i':=\varphi'(p_i')$ .
\qed
\enddemo

Therefore we can prove the rationality of ${\Cal M}_{6,n}$ via standard results
on quotients. To this purpose we consider the incidence variety 
$$
{\Cal E} := \{\ (E,A_1,\dots,A_{n-1},N)\ \vert\  A_{n-1},N \in E\  \}\subseteq \vert \Cal O_\p2(3) \vert
\times (\p2)^{n}
$$
Denoting by $S^9{\Cal E}$ the symmetric product over $ \{\ (E,A_1,\dots,A_{n-1})\ \vert\  A_{n-1}\in E\  \}$,
we notice the existence of a natural chain of rational maps
$\beta\circ\alpha\colon W_{7,9}^n \dasharrow S^9{\Cal E}\dasharrow{\Cal E}$,
where
$$
\alpha(D,A_1,\dots,A_{n-1})=(E,A_1,\dots,A_{n-1},N_1+\ldots+N_9),
$$
$N_1,\dots,N_9$ being the
nodes of $D \in V_{7,9}^{tang}$ on the smooth cubic $E$,
$$
\beta(E,A_1,\dots,A_{n-1},N_1+\ldots+N_9) =
(E,A_1,\dots,A_{n-1},N),
$$
the point $N$ being the unique element in $\left\vert N_1+\ldots+ N_9-8A_1\right\vert$ on the smooth cubic $E$.

\proclaim{Lemma 2.3}
The action of $\PGL_3$ on ${\Cal E}$ has trivial stabilizer at the general point and the quotient ${\Cal E}/\PGL_3$ is rational.
\endproclaim
\demo{Proof}
Let $(E,A_1,\dots,A_{n-1},N)\in{\Cal E}$ be general. Each projectivity leaving fixed both the curve $E$ and the points $A_1,\dots,A_{n-1},N$ induces an automorphism of the abstract elliptic curve $E$ fixing $A_{n-1},N$. Since the automorphism group of the $1$--pointed general elliptic curve $(E,N)$ is isomorphic to ${\Bbb Z}_2$ (see e.g. [Ha], Corollary IV.4.7), then a general choice of $A_{n-1}$ implies that each element in the stabilizer of $(E,A_1,\dots,A_{n-1},N)$ restricts to the identity on $E$, hence has to be the identity on the whole plane since the embedded curve $E\subseteq\p2$ contains $4$ points in general position.

Let $T\subseteq\PGL_3$ be the stabilizer of the points $[1,0,0]$, $[0,0,1]$ and the line $r:=\{\ x_2=0\ \}$ in $\PGL_3$: $T$ is a triangular subgroup of $\PGL_3$.
Clearly
$$
{\Cal E}_0:=\{\ (E,A_1,\dots,A_{n-2},[1,0,0],[0,0,1])\in{\Cal E}\ \vert\ \text{$E$ is tangent to $r$ at $[1,0,0]$}\ \}
$$
is a $(\PGL_3,T)$--section of ${\Cal E}$ in the sense of [Ka] (see also Section 3 of [Do] where it is called \lq\lq slice\rq\rq). In particular ${\Cal E}/\PGL_3\approx{\Cal E}_0/T$ by [Ka], Proposition 1.2.

Define the vector space $V:=\langle x_0^2x_2, x_0x_1^2,x_0x_1x_2,x_0x_2^2,x_1^2x_2,x_1x_2^2,x_2^3\rangle$: then
$$
{\Cal E}_0\cong{\Bbb P}(V)\times(\p2)^{n-2}\cong (V\oplus {\Bbb C}^{\oplus3(n-2)})/({\Bbb C}^*)^{n-1}
$$
thus
$$
{\Cal E}/\PGL_3\approx{\Cal E}_0/T\approx({\Bbb P}(V)\times(\p2)^{n-2})/T\approx(V\oplus {\Bbb C}^{\oplus3(n-2)})/({\Bbb C}^*)^{n-1}\rtimes T.
$$
Since the action of $({\Bbb C}^*)^{n-1}\rtimes T$ on $(V\oplus {\Bbb C}^{\oplus3})$ is linear and triangular, then ${\Cal E}/\PGL_3$ turns out to be rational by [Vi].
\qed
\enddemo

We are now ready to prove the following

\proclaim{Theorem A for $g=6$ and $2\le n\le 8$}
The moduli space ${\Cal M}_{6,n}$ is rational for $2\le n\le 8$.
\endproclaim
\demo{Proof}
Due to the above construction it suffices to prove that $W_{7,9}^n/\PGL_3$ is rational if $2\le n\le 8$. The maps $\alpha,\beta$ are trivially $\PGL_3$--equivariant by construction. Moreover the action of $\PGL_3$ on $S^9 \Cal E$ has trivial stabilizer at the general point since, otherwise,  by projecting down over $\Cal E$ we would obtain points with non--trivial stabilizer.

The typical fibre of $\beta$ over $(E,A_1,A_2,\dots,A_{n-1},N)$ is given by
all divisors $N_1+\ldots+N_9$ such that
$N_1 + \ldots + N_9 \in \left\vert N+8A_1\right\vert$. In particular $\beta$ is dominant and it is a projective bundle over $\Cal E$ with typical fibre $\p{8}$. The section
$\{\ (E,A_1,A_2,\dots,A_n,N+8A_1)\ \} \subset S^9{\Cal E}$ is
$\PGL_3$--equivariant: since $\PGL_3$ is reductive there exists a $\PGL_3$--invariant hyperplane in each fibre, so Lemma 2.3 above and the \lq\lq no--name method\rq\rq\ (see [Do], Section 4) yield the rationality of the quotient $S^9{\Cal E}/ \PGL_3$.

The typical fibre of $\alpha$ over $(E,A_1,A_2,\dots,A_{n-1},N_1+\ldots+N_9)$ is given by
all integral septic curves passing doubly through $N_1,\ldots,N_9$,
simply through $A_1,\dots,A_{n-2}$ and tangent to the conic $E$ at the point $A_{n-1}$ which is a linear system of dimension $8-n$, while the sublocus of non--integral septics has dimension at most $7-n$: it then follows from Bertini Theorem that the general fibre is non-empty and it is isomorphic to $\p{8-n}$. Let $r$ be the tangent line at $E$ in the point $A_{n-1}$. Then the section $\{\ (2E + r,A_1,\dots,A_{n-1})\ \} \subset W_{7,9}^n$ is  $\PGL_3$--equivariant. Thus the same argument as above yields the stated rationality.
\qed
\enddemo

\remark{Comments to Section 2}
If we try to repeat verbatim the above construction
for the case $g=7$, we encounter two problems.
On one hand, the general fibre of $\alpha$ should be given by all
plane octic curves passing doubly through $14$
general points and simply through $4$ general points,
hence it is empty.
On the other hand, the morphism $\beta$ cannot be defined exactly in
the same way since on a smooth quartic curve there
is no group law.
\endremark
\medbreak

\head
3. The unirationality of ${\Cal M}_{8,n}$ for $1\le n\le 11$
\endhead

Let $(C,p_1,\dots,p_{11})\in{\Cal M}_{8,11}$ be general. Due to the Dimension Theorem (see [G--H]: see also [A--C--G--H], (1.5)), the minimum $d$ for which $C$ carries a $g^1_d$ is $5$. Moreover it carries exactly 21 distinct $g^1_5$ by a formula of Castelnuovo (see [Ca]: see also [Fu], Example 14.4.5 or Chapter VII of [A--C--G--H]). By choosing two morphisms $\varphi_i\colon C\to \p1$, $i=1,2$, corresponding to two fixed $g_5^1$'s on $C$, we then obtain a morphism $\varphi:=(\varphi_1,\varphi_2)\colon C\to Q:=\p1\times\p1$.

We first notice that $\varphi$ has either degree $5$ or it is birational onto its image. The first case occurs if and only if $\varphi_1$ and $\varphi_2$ are induced by the same $g_5^1$. Thus we can assume that $\varphi$ is birational onto its image from now on. By construction $\varphi(C)\in\vert\Ofa Q(5,5)\vert$ and it is naturally endowed with an ordered ${11}$--tuple $(A_1,\dots,A_{11})$ of points, namely $A_i:=\varphi(p_i)$. Thus the arithmetic genus of $\varphi(C)$ is $16$, and the genus formula then yields
$$
\sum_{N\in\varphi(C)}{{m_N(m_N-1)}\over2}=16-8=8
$$
where $m_N$ denotes the multiplicity of the point $N\in\varphi(C)$. It follows that $\varphi(C)$ is necessarily singular and it carries at most fourfold points. The linear system of conics through all the points of multiplicity greater than two would cut on $\varphi(C)$ linear series with negative Brill--Noether number. Since $C$ is general this is not possible (again by the Dimension Theorem), thus $\varphi(C)$ carries at most double points as singularities. Let us consider the incidence variety
$$
\align
{\Cal I}:=\{\ &(D,N_1,\dots,N_8,A_1,\dots,A_{11})\in\vert \Ofa Q(5,5)\vert\times Q^{ 19}\ \vert\\
& \text{\ \ $A_i\in D$, $N_j$ is double on $D$}\ \}\mapright\alpha Q^{ 19}
\endalign
$$
(in this notation a point $N_j$ is repeated $t$ times if it is double and there are $t-1$ other double points infinitely near to it). 

\proclaim{Proposition 3.1}
For the general point of $(D,N_1,\dots,N_8,A_1,\dots,A_{11})\in{\Cal I}$, the curve $D$ is irreducible with no singularities but nodes at $N_1,\dots,N_8$, hence its normalization is a ${11}$--pointed curve of genus $8$.
\endproclaim
\demo{Proof}
We have only to check the more or less obvious fact that for the general point $(D,N_1,\dots,N_8,A_1,\dots,A_{11})\in{\Cal I}$ the curve $D$ is irreducible and it has no singularities but the nodes $N_1,\dots,N_8$.

Chose a general $8$--tuple $(N_1,\dots,N_8)\in Q^{ 8}$ of pairwise distinct points. Then there exists a smooth irreducible element $D'\in\vert\Ofa Q(2,3)\vert$ through $N_1,\dots,N_8$. For each other point $A\in Q$ there exists a divisor $D''\in\vert\Ofa Q(3,2)\vert$ such that $D'+D''\in \vert\Ofa Q(5,5)\vert$ is smooth at $A$. Thus the general $D$ has no singularities but $N_1,\dots,N_8$.

Due to the exact sequence 
$$
0\longrightarrow \Ofa Q(1,-1)\longrightarrow \Ofa Q(3,2)\longrightarrow \Ofa Q(1,-1)\otimes\Ofa{D'}\longrightarrow 0
$$
one infers that the very ample linear system $\vert  \Ofa Q(1,-1)\otimes\Ofa{D'}\vert$ is cut out by $\vert  \Ofa Q(1,-1)\vert$. Since the general divisor in $\vert  \Ofa Q(1,-1)\otimes\Ofa{D'}\vert$ has only simple points in its support then we are able to find an element $D''\in \vert\Ofa Q(3,2)\vert$ intersecting transversaly $D'$ at $N_1,\dots,N_8$, whence $D'+D''\in \vert\Ofa Q(5,5)\vert$ has nodes at the points $N_1,\dots,N_8$. Thus the general $D$ has nodes at the points $N_1,\dots,N_8$.

Now an easy computation shows that a curve in $\vert\Ofa Q(5,5)\vert$ with exactly eight nodes as singularities cannot be reducible.
\qed
\enddemo

From the construction above and Proposition 3.1 it follows that

\proclaim{Corollary 3.2}
The composition of the forgetful map onto $\vert \Ofa Q(5,5)\vert\times Q^{ 19}$ with the normalization gives a dominant rational map $m\colon{\Cal I}_n\dashrightarrow{\Cal M}_{19}$.
\qed
\endproclaim

The fibre of $\alpha$ over a point in $Q^{ 19}$ is a possibly empty subspace of $\vert\Ofa Q(5,5)\vert$. Since each double point imposes three linear conditions on the divisors through it, then there exists an open set ${\Cal V}\subseteq Q^{ 19}$ non--intersecting the diagonals of $Q^{ 19}$ over which the fibres of $\alpha$ are projective spaces of constant dimension $\dim\vert\Ofa Q(5,5)\vert-24-11=0$, i.e. a single point: in particular ${\Cal I}$ contains an irreducible component
$$
{\Cal I}^{gen}\subseteq \vert \Ofa Q(5,5)\vert\times Q^{ 19}
$$
of dimension $38$ which is birational via $\alpha$ to a non--empty open subset ${\Cal U}\subseteq{\Cal V}\subseteq Q^{ 19}$, thus it is irreducible and rational.

The fibres of $\alpha$ containing curves of $\vert\Ofa Q(5,5)\vert$ with infinitely near points map over some diagonal in $Q^{ 19}$. Since a double point infinitely near to a double point imposes five linear conditions on the divisors through it then the union of such fibres is a subscheme ${\Cal I}^{sp}\subseteq{\Cal I}$ of dimension at most $37$. 

Since on ${\Cal I}$ there is a natural action of the group of automorphisms $\Aut(Q)$, which has dimension $6$, then the fibres of $m$ have at least dimension $6$, thus ${\Cal I}^{sp}$ maps onto a proper subscheme of ${\Cal M}_{8,11}$. It follows

\proclaim{Theorem B for $g=8$}
The moduli space ${\Cal M}_{8,n}$ is unirational for $1\le n\le 11$.
\endproclaim
\demo{Proof}
The restriction of $m$ to ${\Cal I}^{gen}$ is dominant. Since ${\Cal I}^{gen}$ is rational and dominates ${\Cal M}_{8,11}$, then ${\Cal M}_{8,11}$ turns out to be unirational. Since the forgetful map ${\Cal M}_{8,11}\to{\Cal M}_{8,n}$ is dominant for $n\le10$, it turns out the unirationality of ${\Cal M}_{8,n}$ for $1\le n\le 10$ too.
\qed
\enddemo

\head
4. The unirationality of ${\Cal M}_{10,n}$ for $1\le n\le 3$
\endhead
Let $(C,p_1,p_2,p_3)\in{\Cal M}_{10,3}$ be general: $C$ carries exactly $42$ distinct $g^1_6$ by the formula of Castelnuovo quoted above. In [Se] (see also [A--C])  it is proved the existence  of a morphism $\varphi\colon C\to\p2$ whose image $\varphi(C)\subseteq\p2$ is a curve of degree $9$ with one ordinary triple point and $15$ nodes and no other singularities.
Let us consider the incidence variety
$$
\align
{\Cal I}:=\{\ &(D,N_0,N_1\dots,N_{15},A_1,A_2,A_{3})\in\vert \Ofa{\p2}(9)\vert\times (\p2)^{ 19}\ \vert\\
& \text{\ \ $A_i\in D$, $N_0$ is triple on $D$, $N_j$ is double on $D$ for $j>0$}\ \}\mapright\alpha (\p2)^{ 19}
\endalign
$$
(in this notation a point $N_j$ is repeated $t$ times if it is double and there are $t-1$ other double points infinitely near to it). Corollary~4.7 of [A--C] in the particular case $n=9$, $d=6$ and $\delta=15$ yields the following

\proclaim{Proposition 4.1}
For the general point $(D,N_0,N_1\dots,N_{15},A_1,A_2,A_3)\in{\Cal I}$, the curve $D$ is irreducible with no singularities but an ordinary triple point at $N_0$ and nodes at $N_1,\dots,N_{15}$, hence its normalization is a $n$--pointed curve of genus $10$.
\qed
\endproclaim

From the construction above and Proposition 4.1 it follows that

\proclaim{Corollary 4.2}
The composition of the forgetful map onto $\vert \Ofa{\p2}(9)\vert\times (\p2)^{ 19}$ with the normalization gives a dominant rational map $m\colon{\Cal I}\dashrightarrow{\Cal M}_{10,3}$.
\qed
\endproclaim

As in the genus $8$ case, then ${\Cal I}$ contains an irreducible component
$$
{\Cal I}^{gen}\subseteq \vert \Ofa{\p2}(9)\vert\times (\p2)^{ 19}
$$
of dimension $38$ which is birational to a non--empty open subset ${\Cal U}\subseteq (\p2)^{ 19}$ not intersecting the diagonals of $(\p2)^{ 19}$, thus it is irreducible and rational. Again the locus of points in ${\Cal I}$ carrying infinitely near singularities has dimension at most $37$,
thus

\proclaim{Theorem B for $g=10$}
The moduli space ${\Cal M}_{10,n}$ is unirational for $1\le n\le 3$.
\endproclaim
\demo{Proof}
The restriction of $m$ to ${\Cal I}^{gen}$ is dominant. Since ${\Cal I}^{gen}$ is a projective bundle dominating ${\Cal M}_{10,3}$, then ${\Cal M}_{10,3}$ turns out to be unirational. Since the forgetful map ${\Cal M}_{10,3}\to{\Cal M}_{10,n}$ is dominant for $n\le2$, it turns out the unirationality of ${\Cal M}_{10,n}$ for $1\le n\le 2$ too.
\qed
\enddemo

\remark{Comments to Sections 3 and 4}
Our proof of the unirationality of ${\Cal M}_{g,n}$ depends on the existence of models with a particular configuration of singularities. In the two cases $g=8,10$ we used two different models and it is natural to ask why we made such a choice. 

For instance  we could embed the general curve of genus $g=8$ in $\p2$ as a plane curve of degree $8$ with a triple point and $10$ nodes and prove the unirationality of ${\Cal M}_{8,n}$ as for $g=10$. However, a standard easy parameter count shows that this approach yields the unirationality up to $n=8$ which is in any case greater than Logan's bound but lower than the value obtained with the model on a quadric.

Projecting our model from a double point onto a plane we obtain again a curve of degree $8$ with $2$ triple points and only $7$ double points. Since a triple point imposes on curves of given degree less conditions than $3$ double points, we are able with this model to improve the unirationality result.

With this remark in mind one could hope to construct plane models with more triple points. Unfortunately an easy parameter count shows that the locus in ${\Cal M}_8$ of curves that can be represented as plane curve of degree $8$ with at least $3$ triple points has codimension $1$, so our result is the best possible via our method.

When $g=10$ we can also try to argue as in the case $g=8$ by mapping birationally a general curve $C$ of genus $10$ onto a singular curve of bidegree $(6,6)$ on a smooth quadric. In this case one expects $15$ nodes as singularities but this cannot be proved with the same method we used for $g=8$. Indeed, in this case the Brill--Noether number can be non negative even if triple points occur on the image of $C$. On the other hand, even if we are able to prove that for general $C$ all the singularities are double points and thus even if it makes sense to consider an incidence defined as in the genus $8$ case, it is not evident a priori that the general fibre of such morphism contains irreducible curves of genus $10$. Moreover even if all would work in the right way, such an approach would yield the same bound $n\le3$! 

Thus it is natural to ask if it is possible to represent curves of genus $10$ as  plane curves of degree $9$ with at least two triple points. Again an easy parameter count shows that the locus of such curves in ${\Cal M}_{10}$ has codimension $1$.
\endremark
\medbreak

\head
5. The unirationality of ${\Cal M}_{12,1}$ (and of ${\Cal M}_{14,n}$ for $1\le n\le 2$)
\endhead

We first recall the main steps of Verra's proof of the unirationality of ${\Cal M}_{g}$ for $g=12,14$ (see [Ve]), also in order to fix the notation used throughout this section.

Let $\Hilb_{p(t)}(\p r)$ be the Hilbert scheme parameterizing subschemes $X\subseteq\p r$ with Hilbert polynomial $p(t)$.

\definition{Definition 5.1}
For each triple of positive integers $(d,g,r)$, we denote by ${\Cal C}_{d,g,r}$ the subset of of all smooth curves $C\in\Hilb_{dt+1-g}(\p r)$ such that $\Ofa C(1)$ is non--special and the rank of the natural restriction morphisms $H^0\big(\p r,\Ofa{\p r}(f)\big)\to H^0\big(C,\Ofa C(f)\big)$ is maximal for each $f\in\Bbb Z$.
\enddefinition

Clearly ${\Cal C}_{d,g,r}\subseteq\Hilb_{dt+1-g}(\p r)$ is open and it is non--empty if $3\le r\le d-g$ due to [B--E1], [B--E2], [B--E3]. Moreover, in Section~1 of [Ve], it is also proved the following

\proclaim{Lemma 5.2}
If $g\le 9$ then ${\Cal C}_{d,g,r}$ is irreducible and unirational.
\qed
\endproclaim

Due to the definition of $H^0\big(\p r,\Im_C(f)\big)$, for each $f\in\Bbb Z$ the dimension $h^0\big(\p r,\Im_C(f)\big)$ is constant for $C\in{\Cal C}_{d,g,r}$, thus there exists a positive integer $f$ satisfying $h^0\big(\p r,\Im_C(f)\big)\ge r-1$ for each $C\in {\Cal C}_{d,g,r}$: from now on we will assume that the choice of such an $f$ has been made.

Let ${\Cal G}_{d,g,r}\to{\Cal C}_{d,g,r}$ be the Grassmann bundle parameterizing pairs $(C,V)$ such that $C\in {\Cal C}_{d,g,r}$ and $V\subseteq H^0\big(\p r,\Im_C(f)\big)$ has dimension $r-1$. For each point $(C,V)\in {\Cal G}_{d,g,r}$ let $B_{(C,V)}$ be the base locus of the linear system of hypersurfaces associated to $V$.

\definition{Definition 5.3}
Let $n_{f,d,g,r}:=h^0\big(\p r,\Im_C(f)\big)-(r-1)$ and define the incidence variety
$$
\align
{\Cal I}_{f,d,g,r}:=\{\ &(C,V,A_1,\dots, A_n)\in{\Cal G}_{d,g,r}\times (\p r)^{n_{f,d,g,r}}\ \vert\\
&A_i\in B_{(C,V)}, i=1,\dots,n_{f,d,g,r}\ \}\mapright\alpha{\Cal C}_{d,g,r}\times(\p r)^{n_{f,d,g,r}}.
\endalign
$$
\enddefinition

Points in ${\Cal I}_{f,d,g,r}$ are $(n_{f,d,g,r}+2)$--tuples $(C,V,A_1,\dots, A_{n_{f,d,g,r}})$ such that $C\in {\Cal C}_{d,g,r}$ and $V\subseteq H^0\big(\p r,\Im_C(f)\big)$ has dimension $r-1$. By definition, $C\subseteq B_{(C,V)}$ and $B_{(C,V)}$ could be a curve. If this is the case then $B_{(C,V)}$ is a complete intersection curve and, due to the definition of $n_{f,d,g,r}$, $V$ is exactly the subspace of $H^0\big(\p r,\Im_C(f)\big)$ corresponding to hypersurfaces through $C$ and $A_1,\dots, A_{n_{f,d,g,r}}$. In particular, in this case $\alpha$ is generically injective, thus it is birational since it is surjective by construction.

Moreover, there is also an algebraic linkage $B_{(C,V)}=C\cup D_{(C,V)}$, where $D_{(C,V)}$ is also a curve. If $B_{(C,V)}$ is nodal and $D_{(C,V)}$ is smooth, irreducible and non--degenerate in $\p r$ then its degree $d'$ and genus $g'$ are given by the formulas $d+d'=f^{r-1}$ and $2(g-g')=(f(r-1)-r-1)(d-d')$ (see [Mi], Corollaries~5.2.13 and 5.2.14).

Let us restrict ourselves to the case $(f,d,g,r)=(2,15,9,6)$. Notice that in this case $h^0\big(\p 6,\Ofa{\p 6}(2)\big)=28>22=h^0\big(C,\Ofa C(2)\big)$: since $C$ has maximal rank  then
$$
n_{2,15,9,6}:=h^0\big(\p 6,\Im_C(2)\big)-5=h^0\big(\p 6,\Ofa{\p 6}(2)\big)-h^0\big(C,\Ofa C(2)\big)-5=1.
$$
\proclaim{Proposition 5.4}
There exists an open subset ${\Cal U}\subseteq{\Cal C}_{15,9,6}\times\p 6$ such that for each $(C,V,A_1)\in {\Cal I}^{gen}_{2,15,9,6}:=\alpha^{-1}({\Cal U})$, then $D_{(C,V)}$ is smooth, irreducible and non--degenerate in $\p 6$ with degree $17$ and genus $12$.
\endproclaim
\demo{Proof}
There is a dominant morphism $g\colon{\Cal I}_{2,15,9,6}\to{\Cal G}_{15,9,6}$. In Section 7 of [Ve] it has been proved the existence of an open set ${\Cal V}\subseteq{\Cal G}_{15,9,6}$ such that for each $(C,V)\in {\Cal V}$ the curve $D_{(C,V)}$ satisfies the assertion, hence the same is true for each general point $(C,V,A_1)\in g^{-1}({\Cal V})$.
\qed
\enddemo

In particular, as in the previous sections, we infer that ${\Cal I}^{gen}_{2,15,9,6}$ is birational via $\alpha$ to ${\Cal U}$ which is unirational due to Lemma~5.2. We are now ready to prove our

\proclaim{Theorem B for $g=12$}
The moduli space ${\Cal M}_{12,1}$ is unirational.
\endproclaim
\demo{Proof}
Thnaks to Proposition~5.4, we have a rational map $m_{2,15,9,6}\colon{\Cal I}^{gen}_{2,15,9,6}\to{\Cal M}_{12,1}$, hence it suffices to check it is dominant. 

Let $(D,p)\in {\Cal M}_{12,1}$ be general. Then $D\in {\Cal M}_{12}$ is general too, thus [Ve], Theorem~4.5 (see Section~7 for its application in this case) implies the existence of $(C,V)\in {\Cal G}_{15,9,6}$ such that $D_{(C,V)}= D$. Let $A\in D_{(C,V)}\subseteq\p6$ be the point corresponding to $p$ via such an isomorphism: since $A\in D_{(C,V)}\subseteq B_{(C,V)}$ and it is general on $D_{(C,V)}$, then $p\not\in C$, $V$ is exactly the space associated to the linear system of hypersurfaces of degree $2$ through $C$ and $A$, i.e. $(C,V,A)\in {\Cal I}_{2,15,9,6}^{gen}$. By construction $m_{2,15,9,6}(C,V,A)=(D,p)$.
\qed
\enddemo

\remark{Remark 5.5}
In [B--V], Theorem~4.7 the authors prove with the above procedure the unirationality of ${\Cal M}_{14,n}$ for $1\le n\le2$. Indeed, due to the result of Section 6 of [Ve], Proposition~5.4 holds true also in the case $(f,d,g,r)=(2,14,8,6)$. In this case $n_{2,14,8,6}=2$ and $D_{(C,V)}$ is a curve of genus 14, thus one can define the rational map $m_{2,14,8,6}\colon{\Cal I}^{gen}_{2,14,8,6}\to{\Cal M}_{14,2}$. As in the case $g=12$ one then checks that $m_{2,14,8,6}$ is dominant.  
\endremark
\medbreak

\remark{Comments to Section 5}
One could hope to generalize the above construction also for other genera. Reading [Ve] one checks that this is possible if the hypotheses of Theorem 4.5 of [Ve] are satisfied. This occurs if the sequence of  non--negative integers $(r,f,d,g,d',g')$ satisfies the following system 
$$
\cases
f={{r+2}\over{r-2}}\in{\Bbb Z}\\
d+d'=f^{r-1}\\
d-g=r-1\\
2(g-g')=((r-1)f-r-1)(d-d')\\
g'-(r+1)(g'-d'+r)\ge0.
\endcases
$$
Solving the above system in $\Bbb Z$ one then finds all its solutions, which are
$$
\gather
(6,2,15,9,17,12),\quad(6,2,14,8,18,14),\quad(6,2,16,10,16,10),\\(4,3,12,8,15,14),\quad(4,3,13,9,14,11),\\(3,5,12,9,13,12).
\endgather
$$
The first two ones have been used in the present section to prove the unirationality of ${\Cal M}_{12,1}$ and ${\Cal M}_{14,2}$. In the third case $n_{2,16,10,6}=0$, we cannot apply the above method with the fourth sequence and the result in Section 4 remains the best possible with the method used in this paper.

Since, we already know (see the table with the values of the functions $\sigma$ and $\tau$ in the introduction) that ${\Cal M}_{11,10}$ is unirational but $\kappa({\Cal M}_{11,10})\ge0$, we cannot improve the known result in genus $g=11$ applying the above method in the fourth case $(r,f,d,g,d',g')=(4,3,13,9,14,11)$.

The fifth sequence has been used in [Ve] to prove the unirationality of ${\Cal M}_{13}$: unfortunately the method of the proof in [Ve] does not allow us to fix points in order to infer the unirationality ${\Cal M}_{13,n}$ for whatever value of $n$. 

Consider finally the last sequence. Since $n_{f,d,g,r}=2$ in this case, the above method would lead to the unirationality of ${\Cal M}_{12,2}$, improving our result. Therefore, it would be very interesting to understand if at least one of the sufficient conditions of Theorem~4.5 in [Ve] is satisfied.

The first condition is not satisfied since the homogeneous ideal $\oplus_{n\in\Bbb Z}H^0\big (\p3,\Im_{C}(n)\big)$ in $\p3$ of a general point $C$ of the family ${\Cal C}_{12,9,3}$ is minimally generated by sextics other than quintics. Indeed $h^0\big(C,\Ofa C(5)\big)=52$, $h^0\big(C,\Ofa C(6)\big)=64$, $h^0\big(\p3,\Ofa{\p3}(5)\big)=56$,  $h^0\big(\p3,\Ofa{\p3}(6)\big)=84$, hence the maximal rank condition yields $h^0\big(\p3,\Im_{C}{\p3}(5)\big)=4$,  $h^0\big(\p3,\Im_{C}(6)\big)=20$.

Thus the only chance is to verify the second condition. A careful reading of Verra's construction will convince the interested reader that instead of checking the irreducibility of the universal Brill--Noether locus ${\Cal W}_{13,12}^3$, it suffices to use the irreducibility of the subscheme ${\Cal J}'(13,12,3)$ of the Hilbert scheme of curves in $\p3$, whose general point corresponds to a smooth, irreducible and non--degenerate curve of degree $13$ and genus $12$ (see [K--K], Theorem~2.7~(iii)).

Hence the crucial point it seems to establish if ${\Cal G}_{12,9,3}$ satisfies the key condition. In principle one could try to apply [Ve], Proposition~3.4, but the assumption that the twisted sheaf $\Im_{C}(5)$ of ideals of a general point $C$ in ${\Cal C}_{12,9,3}$ is globally generated is hard to be verified. Another possibility could be to address directly the key condition but again this seems to be completely out of reach for us. Therefore we leave it as an intriguing open problem.
\endremark
\medbreak

\enddocument